\documentclass[12pt]{article}
\usepackage{fullpage, amsthm}
\usepackage{hyperref}
\usepackage{graphicx}
\usepackage{amssymb}
\usepackage{epstopdf}
\usepackage{amsmath}

\DeclareSymbolFont{symbolsC}{U}{pxsyc}{m}{n}
\DeclareMathSymbol{\coloneqq}{\mathrel}{symbolsC}{"42}

\newcommand{\fracpar}[2]{\frac{\partial #1}{\partial #2}}
\newcommand{\be}{\begin{equation}}
\newcommand{\ee}{\end{equation}}
\newcommand{\bea}{\begin{eqnarray}}
\newcommand{\eea}{\end{eqnarray}}
\newcommand{\eps}{\epsilon}

\newcommand{\R}{\mathbb{R}}

\newcommand{\x}{{x}}
\newcommand{\J}{{J}}
\newcommand{\Haus}{\mathcal{H}}
\newcommand{\bsigma}{{\sigma}}
\newcommand{\ubar}{\bar{u}}

\newtheorem{theorem}{Theorem}[section]
\newtheorem{lemma}[theorem]{Lemma}

\newtheorem{proposition}[theorem]{Proposition}
\newtheorem{corollary}[theorem]{Corollary}
\newtheorem{definition}[theorem]{Definition}

\newenvironment{remark}[1][Remark]{\begin{trivlist}
\item[\hskip \labelsep {\bfseries #1}]}{\end{trivlist}}
\title{Current Density Impedance Imaging of  an Anisotropic Conductivity in a Known Conformal Class}
\author{  
{Nicholas Hoell\footnote{Department of Mathematics, University of Toronto, Toronto, Ontario, Canada M5S
2E4. E-mail: nhoell@math.toronto.edu.  The author is  supported in part by a MITACS Postdoctoral Fellowship.} }\qquad
{Amir Moradifam\footnote{Department of Applied Physics and Applied Mathematics, Columbia University, New York, NY, USA. E-mail: am3937@columbia.edu. The author is  supported in part by  Mitacs and NSERC Postdoctoral Fellowships. }
\qquad Adrian Nachman\footnote{Department of Mathematics and the
Edward S. Rogers Sr. Department of Electrical and Computer
Engineering, University of Toronto, Toronto, Ontario, Canada. E-mail: nachman@math.toronto.edu. The author is supported in part by an NSERC Discovery Grant.}\qquad }}
\date{February 28, 2013}

\begin{document}
\maketitle

\begin{abstract} {We present a procedure for recovering the conformal factor of an anisotropic conductivity matrix in a known conformal class in a domain in $\R^n$ with $n\geq 2$.    The method requires one internal measurement, together with a priori knowledge of the conformal class (local orientation) of the conductivity matrix.  This problem arises in the coupled-physics medical imaging modality of Current Density Impedance Imaging (CDII) and the assumptions on the data are suitable for measurements determinable from cross-property based couplings of the two imaging modalities CDII and Diffusion Tensor Imaging (DTI).   We show that the corresponding electric potential is the unique solution of a constrained minimization problem with respect to a weighted total variation functional defined in terms of   the physical data.  Further, we show that the associated equipotential surfaces are area minimizing with respect to a Riemannian metric obtained from the data. The results  are also extended to allow the presence of perfectly conducting and/or insulating inclusions.  }
\end{abstract}

\section{Introduction}
Biological tissues such as muscle or nerve fibres are known to be electrically anisotropic (see e.g. \cite{nicholson, roth}). In this paper, we consider the problem of recovering an anisotropic electric conductivity $\sigma$ of a body $\Omega$ from measurement of one current $J$ in the interior. Such interior data can be obtained by Current Density Imaging (CDI), a method pioneered at the University of Toronto (\cite{joy1, joy3}) that used a Magnetic Resonance Imager (MRI) in a novel way. We also rely on the MRI-based Diffusion Tensor Imaging (DTI) method to determine the conformal class of  $\sigma$, as in the new DT-CD-II method recently introduced and tested experimentally in (\cite{weijing}). Thus, we assume that the matrix-valued conductivity function is of the form:
\begin{equation}\label{sig2}
\sigma(\x)=c(\x)\sigma_0(\x),
\end{equation}
 with $\sigma_0(\x)$  known  from DTI, and with the cross-property factor $ c(\x)$ a scalar function to be determined. This assumption is motivated  by a number of physical studies which have shown a linear relationship between the conductivity tensor and the diffusion tensor (see e.g. \cite{tuch, weijing} and further references therein). 
 
We show that, in dimension $n \geq 2$, the cross-property factor $ c(\x)$ can be determined from knowedge of the current $J$ in $\Omega$  and of the corresponding prescribed voltage $f$ on the boundary $\partial \Omega$. In fact, the only internal data we require is the scalar function 
\begin{equation}\label{sca}
a=(\sigma_0^{-1} J \cdot J)^{\frac{1}{2}}
\end{equation}
(with $ \sigma_0^{-1} $ denoting the inverse of the matrix $\sigma_0$). This turns out to be the appropriate extension of the corresponding earlier result for isotropic conductivities (\cite{alex3}), where the interior data was the magnitude $|\J|$.

The method we will be presenting is based on the minimization of a weighted total variation functional defined in terms of $a(\x)$and $\sigma_0(\x)$. (See Theorem \ref {mainthm} for the precise statement).

More generally, we will show that when $\Omega$ contains perfectly conducting and/or insulating inclusions, then knowledge of $a$,  $\sigma_0$ and $f$ determines the location of these inclusions (except in exceptional cases), as well as the function $ c(\x)$  (hence also the anisotropic conductivity $\sigma$) in their complement.

\subsection{Background and Motivation}
Mathematical work on non-invasive determination of internal conductivity has focused largely on the classical method of Electrical Impedance Tomography (EIT).  There have been major advances in the understanding of this nonlinear inverse boundary value problem (see \cite{uhlmann12} for an excellent review; in particular, see \cite {dksu, ksu} for recent results on recovering anisotropic conductivities in a given conformal class for the special case of admissible manifolds).  It has also been shown that the EIT problem is exponentially ill-posed, yielding images of low resolution (see \cite{isaacson, mandache}).

In a new class of inverse problems (that includes the one studied here) one seeks to overcome the limitations of the reconstructions obtainable from classical boundary measurements by using data that can be measured noninvasively in the  interior of the object. These are known in the literature as 'hybrid problems' (also as 'coupled physics',  'interior data' or 'multi-wave' problems), as they couple two imaging modalities to obtain internal measurements. For overviews of such methods see \cite{bal3, kuchment}. For imaging the electric conductivity, there are several approaches that combine aspects of EIT with MRI: MREIT, CDII, Electric Properties Imaging (see \cite{alex1, seo} for recent reviews) or with ultrasound measurements: Acousto-Electrical Tomography (\cite{wang, abctf, kk}), Impedance-Acoustic Tomography (\cite{scherzer}). 

The starting point for the method presented here is the measurement of one applied current $J(\x)$ at all points $x$ inside a region $\Omega$. We briefly recall the influential idea of \cite{joy1, joy3} for obtaining such interior measurements using MRI. The current $J$ induces a magnetic field  ${B(\x)}$. The component of $B$ parallel to the static field of the imager can be determined at any point inside $\Omega$ from the corresponding change in the phase of the measured magnetization at that location. By performing rotations of the object and repeating the experiment with the same applied current, all three components of $B$ can be recovered, and ${J}(\x)$ is then computed using Amp\'ere's law:
\[\J(\x)=\frac{1}{\mu_0}\nabla \times {B}(\x) \]
where $\mu_0$ is the magnetic permeability (essentially constant in tissue).  
For our purposes, it is important to note that this Current Density Imaging (CDI) method works equally well in anisotropic media, as no knowledge of the conductivity is needed for the determination of the current density $J(\x)$.\\

\label{background}
Inside the body being imaged the electric potential $u(\x)$ corresponding to the voltage  $f(\x)$ on the boundary solves the following Dirichlet problem for the conductivity equation: 
\begin{align}
\label{cond2}
\nabla \cdot \sigma \nabla u&=0, \qquad \x \in \Omega \subset \R^n \\
u\left. \right|_{\partial \Omega}&=f \nonumber
\end{align} 
where $\sigma$ is the (generally tensorial) conductivity of the material.  In the case of isotropic conductivities, (i.e. scalar $\sigma$) considered in \cite{alex1, alex2, alex3, alex4, alex5} and in the absence of insulating or perfectly conducting inclusions one can replace $\sigma$ in the above equation using Ohm's law $|\J|=\sigma |\nabla u|$ to obtain the quasilinear, degenerate, elliptic, variable coefficient $1$-Laplacian equation:
\begin{align}
\label{1Lap}
\nabla \cdot (|\J| \frac{\nabla u}{|\nabla u|})=0, \qquad \x \in \Omega.
\end{align} 
%\be
%\nabla \cdot (|\J| \frac{\nabla u}{|\nabla u|})=0, \qquad \x \in \Omega.
%\ee
The above equation was first introduced, with the above derivation, in the article \cite{kim}, where the Neumann problem was considered and examples of non-existence and non-uniqueness were given to explain that additional data was needed for determining the conductivity. In the article \cite{alex4} it was shown that equipotential surfaces, namely the level sets of $u(x)$, are minimal surfaces with respect to the conformal metric $|\J|^{\frac{2}{n-1}} I_n$, with $I_n$ the $n \times n$ identity matrix; this result was then used to treat the Cauchy problem for equation \eqref{1Lap}. The Dirichlet problem for equation \eqref{1Lap} can also have infinitely many solutions (see \cite{alex3}). A solution around this difficulty was found in \cite{alex3}, where the approach via the partial differential equation \eqref{1Lap} was replaced by the study of the variational  problem for which it is the Euler-Lagrange equation. It was shown that the solution of \ref{cond2} is the unique minimizer for this problem. We recall these results in the following theorem.

\begin{theorem}(\cite{alex3})\label{tim}{Let $\Omega \subset \R^n$, $n \geq 2$ be a domain with a connected $C^{1,\alpha}$ boundary, $\alpha >0$, and let $\mu$ denote Lebesgue measure on $\Omega$ .  Let $(f,|\J|) \in C^{1,\alpha}(\partial \bar{\Omega}) \times C^\alpha(\Omega)$ with $|\J| \neq 0$ $\mu$-a.e.  be associated with an unknown conductivity $\sigma \in C^\alpha(\bar{\Omega})$. 
Then 
\[u_\sigma = \underset{v \in W^{1,1}(\Omega)\bigcap C(\bar{\Omega})}{\text{argmin}} \{ \int_\Omega |\J| |\nabla v|\mu(dx);  \quad v \left.\right|_{\partial \Omega} =f \text{ and } \mu(\{\nabla v = 0\})=0 \} \]
exists and is unique.  

Furthermore,  $\sigma = \frac{|\J|}{|\nabla u_\sigma|} \in L^\infty(\Omega)$ is the unique $C^\alpha(\bar{\Omega})$ scalar conductivity associated to the pair $(f, |\J|)$. }
\end{theorem}
 A generalization of the above result was later obtained in the article \cite{alex5}, where the isotropic conductivity was shown to be determined from knowledge of $|J|$  on the complement  of open regions on which $\sigma$ may be zero (in the case of insulating inclusions) or infinite (for perfectly conducting inclusions).  We refer the interested reader to details in \cite{alex5}.

\subsection{ Statement of Results}
 
In this article we will extend the imaging method described above to the case in which the conductivity is anisotropic and known to be of the form $\sigma(\x)=c(\x)\sigma_0(x)$ where $c(\x)$ is an unknown scalar function and $\bsigma_0 \in C^0(\Omega, Mat(\R,n))$ is a symmetric positive definite matrix-valued anisotropic term, assumed known.  

We shall first prove an anisotropic analogue to Theorem \ref{tim} as a prelude to the more general results accounting for inclusions and less restrictive function spaces.  For this, we will need to precisely    define the class of data that arises from physical measurements. 
\begin{definition}[First notion of admissibility]\label{admis}{A triplet $(f,\sigma_0, a) \in H^{\frac{1}{2}}(\partial \Omega) \times L^\infty_+(\Omega,Mat(\R,n))\times L^2(\Omega)$ shall be said to be \textit{admissible} if there exists a $c(x) \in L^\infty_+(\Omega)$  such that\[a=(\sigma_0^{-1} J \cdot J)^{\frac{1}{2}},\]
where the current
\[J=-c\sigma_0 \nabla u\]
corresponds to  the potential $u\in H^1(\Omega)$, the weak solution to the BVP 

\begin{equation}
\label{bvp1}
\left\{ \begin{array}{ll}  \nabla \cdot (c\sigma_0 \nabla u)=0, & x \in \Omega \\  u\left.\right|_{\partial \Omega} =f. &  \end{array}
\right.
\end{equation}
} 
\end{definition}
The first theorem we prove concerns the minimization of the following functional: 
\be
\label{funcy}
\mathcal{F}[v] \coloneqq \int_{\Omega} a(\sigma_0 \nabla v \cdot \nabla v)^{\frac{1}{2}}d\mu
\ee
where $\mu$ denotes the Lebesgue measure on $\Omega$. We present the following uniqueness result in section \ref{noinclusions}.

\begin{theorem}\label{firstthm}{Let $\Omega \subset \R^n$ be a domain  with a connected $C^{1,\alpha}$ boundary, $\alpha >0$,  and let the triplet $(f,\sigma_0, a) \in C^{1, \alpha}(\partial \Omega)\times C^\alpha(\overline{\Omega},Mat(\R,n))\times C^\alpha(\Omega)$ be an admissible triplet as in Definition \ref{admis} with $|J|>0$ Lebesgue $-a.e.$ in $\Omega$. Denote by $\sigma = c\sigma_0 \in C^\alpha(\overline{\Omega}, Mat(\R, n))$ the unknown generating conductivity for this triplet. 
Then the following minimization problem 
\be
\label{argmin}
\underset{v\left.\right|_{\partial \Omega}=f  }{argmin } \ \{  \mathcal{F}[v];\quad  v \in W^{1,1} \cap C(\overline{\Omega}),\ \mu\{ \nabla v=0\}=0 \}
\ee
has a unique solution $u_\sigma$, where the functional $\mathcal{F}$ is as in \eqref{funcy}. 

Furthermore, the unique $C^{\alpha}(\overline{\Omega}, Mat(\R, n))$ conductivity generating the current density $J$ while maintaining the boundary voltage $f$ is given by $\sigma=c(x)\sigma_0(x)$ with the conformal factor $c$ determined from the  formula
\[c=\frac{a}{(\sigma_0 \nabla u_\sigma \cdot u_\sigma)^{\frac{1}{2}}}.\]
}
\end{theorem}

Following this we establish, in the remainder of section \ref{noinclusions}, the geometrical result that equipotential sets $u^{-1}(\lambda)\coloneqq \{x; u(x)=\lambda \}\cap \overline{\Omega}$ are in fact minimal surfaces with respect to a certain Riemannian metric on $\Omega$ which is defined in terms of  of $\sigma_0(x)$ and $a(x)$.

After  the above preliminary results, we prove the following main uniqueness result in section \ref{inclusions}, which allows for inclusions and more general function spaces.  The precise statement requires an appropriate notion of admissibility, formally defined in section \ref{inclusions} which involves some technical extensions of the criteria in Definition \ref{admis}.  It also requires an extension of the functional $\eqref{funcy}$ appearing in Theorem \ref{firstthm}. A full definition of the proper generalization is postponed until section \ref{prelim} but it takes the form of a weighted total variation $\int_{\Omega}|Dv|_\phi$, where the weight $\phi$ is  defined in terms of $a$ and  $\sigma_0$ and where $|Dv|_\phi$ is a weighted distributional gradient discussed in section \ref{prelim}. 

Our uniqueness result also requires certain natural assumptions on the regions of perfect and zero conductance $O_\infty$ and $O_0$, respectively, as is discussed in greater detail in that section. Further, the zero set of the function $a$, denoted $S$, are also assumed to have a certain topological feature discussed in section \ref{prelim}.  
 
\begin{theorem} \label{unique}
Let $\Omega \subset \R^n$, $n\geq 2$, be a domain with connected
Lipschitz boundary and let $(f, \sigma_0, a)\in H^{1/2}(\partial \Omega)\times C^{0}(\Omega \setminus (O_\infty\cup O_0), Mat(\R,n)) \times
L^{\infty}(\Omega)$ be an admissible triplet generated by an
unknown  conductivity $\sigma$ in the sense of Definition \ref{def}. Then the potential $u$ is a minimizer of the problem
\begin{equation}\label{min_prob}
u=\hbox{argmin} \{\int_{\Omega} |Dv|_{\phi}: v \in BV_{loc}(\Omega \setminus \bar{S}) \ \ \hbox{and} 
\ \ v|_{\partial \Omega}=f\},
\end{equation}
and if $\bar{u}$ is another minimizer of the above problem, then
$\bar{u}=u$ in $\Omega \backslash \{|\J|=0\}$.
 Consequently \[\sigma=\frac{a}{(\sigma_0 \nabla u \cdot \nabla u) ^{\frac{1}{2}}} \sigma_0 \in
C^{\alpha}(\Omega\setminus \overline {Z})\] is the unique conductivity generating the admissible data triplet $(f, \sigma_0, a)$ and $Z$ is an open set discussed later which accounts for the inclusions. 
\end{theorem}

With this shown, we lastly prove  that level sets of such solutions as above minimize the area functional
\[\mathcal{A}(\Sigma)=\int_{\Sigma}a dS \]
over a suitable space of variations of subsets of  $\Sigma \subset \Omega$.

An outline of this paper is as follows.  In section \ref{noinclusions} we present the main ideas building up to Theorem \ref{mainthm}.  Beginning with a formal statement of the problem, we describe  the data of interest in our first definition of admissibility. We then  proceed to establish the uniqueness result under this assumption.  This section avoids  the more technical details that arise from consideration of singular, perfectly conductive, and/or insulating inclusions and also considers potentials satisfying a more restrictive amount of smoothness.  These more restrictive  assumptions help clarify the main  ideas that will be useful later.  In section \ref{prelim} we will introduce some tools from geometric measure theory needed to define an appropriate  weighted space of bounded variation for our minimization problem.  Most of the results presented in this section originated in the article \cite{AB}. In section \ref{inclusions} we formulate a more general notion of admissibility and present the proof of our main uniqueness result, Theorem \ref{unique}.  We then prove, in section \ref{geometric}, that equipotential sets  minimize an area functional defined in terms of  current density measurements.  Finally, some of the technical facts on existence and uniqueness of solutions to a limiting form of the  conductivity equation, suitable in the presence of  inclusions with zero or infinite conductivity, as well as an equivalent optimization result in this setting, are briefly presented  in section \ref{appendix}.  Section \ref{conc}  presents conclusions and acknowledgments.

\section{Anisotropic Current Density Impedance Imaging}
\label{noinclusions}
In this section we present a simplified exposition of the main results of this paper , in order to illustrate the basic ideas used in the argument and to motivate the more general results to be presented later.  We also use this section to briefly introduce some of the key geometric measure-theoretic concepts we will need and expand upon later; some excellent references thereon   may be found in \cite{evans2, federer, giusti, maggi, morgan}.

\subsection{ Uniqueness in the Variational Problem for Inversion}
Assume that 
the conductivity $\sigma$ is of the form $c(\x) \sigma_0(x)$ with $c(x), (\sigma_0)_{ij}(x) \in C^\alpha(\overline{\Omega})$, $\alpha >0$,  $c(x)>0$  and $\sigma_0$ symmetric and positive-definite throughout $\Omega$.

   Throughout the paper we will be using the notation
\be
(\xi, \eta)_{\sigma_0} \coloneqq (\sigma_0 \xi)\cdot \eta,   \qquad |\xi|_{\sigma_0}\coloneqq ((\sigma_0 \xi)\cdot \xi)^{\frac{1}{2}}, \qquad \quad \xi, \eta \in C^{\alpha}(\bar{\Omega}, \R^n)
\ee
to denote the inner product induced by  $\sigma_0$, and the corresponding norm, where $\cdot$ will always be taken to denote the Euclidean dot product.  We also define the space $W^{1,1}_+(\Omega)$ as 
\[W_+^{1,1}(\Omega) \coloneqq \{v \in W^{1,1}(\Omega), \quad \mu(\{x, \nabla v=0 \})=0\}\]   In what follows $\nabla$ denotes the usual (i.e. non-covariant) partial differentiation and we use the Einstein summation convention over repeated upper/lower indices. 
 
  We begin by showing that the solution $u$ to the BVP (\ref{bvp1}) is a minimizer of an action on $\Omega$ that is defined in terms of the internal density magnitude $|J|_{\sigma_0^{-1}}$.  This generalizes the corresponding result for isotropic conductivities in \cite{alex3}.

\begin{lemma}\label{actionlemma}{Assume that $(f,\sigma_0,a)$ is an admissible triplet and let $u$ be a solution to the corresponding forward problem (\ref{bvp1}). Then $u$ is a minimizer of the action integral $\mathcal{F}[\ \cdot \ ]$ defined by the following 
\be
\label{action}
\mathcal{F}[v] \coloneqq \int_\Omega a(\x)|\nabla v|_{\sigma_0}   d\x, \ee
i.e. the relation \be
\mathcal{F}[v] \geq \mathcal{F}[{u}]
\ee
Holds for  all $ v \in W^{1,1}_+(\Omega)$ satisfying $v\left.\right|_{\partial \Omega}=f$.}

\begin{proof}
Let $v \in W^{1,1}_+(\Omega)$.   Since $a$ comes from an admissible triplet, there is a choice of $c(\x)$ such that $a(\x)$ takes the form $a=|\J|_{\sigma_0^{-1}}$.  Then
\begin{align}
\mathcal{F}[v]&= \int_\Omega |\J|_{\sigma_0^{-1}} |\nabla v|_{\sigma_0} d\mu\nonumber \\
&=\int_\Omega c(\x) |\nabla {u}|_{\sigma_0} |\nabla v|_{\sigma_0}  d\mu\nonumber\\
&\geq \int_\Omega c(\x) |(\nabla {u}, \nabla v)_{\sigma_0}|d\mu \label{cauchy}\\
&= \int_\Omega \bsigma \nabla {u} \cdot  \nabla v d\mu \nonumber \\
&=\label{inty} \int_{\partial \Omega}f \bsigma \fracpar{{u}}{\pmb{n}} d{S}\\
&=-\int_{\partial \Omega} f \J \cdot \pmb{n}d{S} \nonumber 
\end{align}
with $\pmb{n}$ an outer-oriented normal to $\partial \Omega$ and where, in line (\ref{inty}), we have integrated by parts and applied the conductivity equation on ${u}$. Equality holds in line (\ref{cauchy}) if and only if $\nabla {u}$ and $\nabla v$ are parallel $\mu -a.e$.   In particular, we have
\begin{align}
\mathcal{F}[{u}]=-\int_{\partial \Omega} f\J \cdot \pmb{n}d{S} \nonumber  
\end{align}
which, on comparing with the above, shows that $u$ is a minimizer, as claimed.

\end{proof}
\end{lemma}

In order to prove the main result of this section we shall need to recall some basic notions from geometric measure theory. Firstly,  by $\mathcal{H}^{d}(\Sigma)$ we denote the $d$-dimensional Hausdorff measure of a set $\Sigma \subset \Omega$ defined as 
\[\Haus^d(\Sigma) \coloneqq \lim_{\delta \downarrow 0 }\text{inf} \{ \sum_{j=1}^n (\text{diam} E_j)^d, \quad \bigcup_{j \in \mathbb{N}} E_j \supset \Sigma, \quad \text{diam} E_j \leq \delta \}  \]
The super-level set of a non-negative function $u(x) \in W^{1,1}(\Omega)$, given by $E_t \coloneqq \Omega \cap \{u>t \}$ has so-called \textit{locally finite perimeter}, in the sense that the vector-valued Radon measure $\nabla \chi_{E_t}$ satisfies $|\nabla \chi_{E_t}|<\infty$ for almost all $t$. For such sets we shall be concerned with the \textit{reduced boundary}.
\begin{definition}{The  \textit{reduced boundary} $\partial^*E$ of a set with locally finite perimeter is the set of points in $\R^n$ for which the following hold;  
 \begin{enumerate}
 \item For all $\eps>0$ one has $\int_{B(x,\eps)}|\nabla \chi_E|>0$
 \item The measure-theoretic outer normal $\nu(x)$ determined by \[ v(x) \coloneqq - \lim_{\eps \downarrow 0}\frac{\int_{B(x,\eps)}\nabla \chi_E}{\int_{B(x,\eps)}|\nabla \chi_E|}\]
 exists, and satisfies $|\nu(x)|=1$.
 \end{enumerate}
}\end{definition}
For a super-level set $E_t$ the unit normal $\nu_t(x)$ exists $\Haus^{n-1}-a.e$ $x \in \partial^* E_t$ (see the remarks in \cite{alex3}).  

We now present the main result of this section.
\begin{theorem}{\label{mainthm}}{Let $\Omega \subset \R^n$, $n\geq 2$, be a domain  with a connected $C^{1,\alpha}$ boundary, $\alpha >0$, and let $(f,\sigma_0, |J|_{\sigma_0^{-1}}) \in C^{1, \alpha}(\partial \Omega)\times C^\alpha(\overline{\Omega},Mat(\R,n))\times C^\alpha(\Omega)$ be an admissible triple with $|J|>0$ $\mu -a.e.$ in $\Omega$. Denote by $\sigma \in C^\alpha(\overline{\Omega})$ the unknown generating conductivity for this triplet. 

 Then the following minimization problem 
\be
\label{argmin}
\underset{v\left.\right|_{\partial \Omega}=f  }{argmin } \{  \mathcal{F}[v]; v \in W^{1,1}_+\cap C(\overline{\Omega})\}
\ee
has a unique solution $u_\sigma$.  

Further, the unique $C^{\alpha}(\overline{\Omega})$ conductivity generating the local current density $J$ while maintaing a boundary voltage $f$ is given by $\sigma=c(x)\sigma_0(x)$ with the conformal factor $c$ determined by
\[c(x)=\frac{|J|_{\sigma_0^{-1}}}{|\nabla u_\sigma|_{\sigma_0}}.\]
}
\begin{proof}
{The proof is essentially the same as the proof of Theorem \ref{tim} given in \cite{alex3} so we present a self-contained but tersely abbreviated proof and refer the reader to that article for full details. Note first that since the triplet is assumed admissible the previous lemma ensures that ${argmin}_{u\left.\right|_{\partial \Omega}=f} \{\mathcal{F}[u]; u \in W^{1,1}_+\cap C(\overline{\Omega})\}$ is nonempty. We choose a minimizer and call it $u_\sigma$.    

To show  uniqueness, assume to the contrary that another minimizer to problem \ref{argmin}, say $\tilde{u} \in W_+^{1,1}\cap C(\overline{\Omega})$, exists.  Recalling  the proof of Lemma \ref{actionlemma} one sees the Cauchy inequality used in \eqref{cauchy} ensures that $\nabla u_\sigma=\lambda(x) \nabla \tilde{u}$ for some non-negative $\lambda$ $\mu-a.e.$.  From this it follows that
\[\frac{\nabla u_\sigma}{|\nabla u_\sigma|}=\frac{\nabla \tilde{u}}{|\nabla \tilde{u}|} \]
holds $\mu-a.e.$  We  show that this  implies equality of the minimizers away from Lebesgue-negligible sets. 

In view of Lemma 2.2 of \cite{alex3},  the above yields an identification of measure-theoretic normal vectors a.e.. It follows that the super-level set $E_t = \{\tilde{u}>t\} \cap \Omega$ has a measure-theoretic normal $\nu_t(x) = -\frac{\nabla \tilde{u}}{|\nabla \tilde{u}|}$ which is continuously extendible from the reduced boundary $\partial^*E_t \cap \Omega$ to the topological boundary $\partial E_t \cap \Omega$. It then follows from a result of De Giorgi (e.g. 4.11 in \cite{giusti}) that, for almost all $t$, the region $\partial E_t \cap \Omega$ is a $C^1$-hypersurface with unit normal $\nu_t(x)$.  From this, a parameterization of each connected component of this hypersurface shows  that $u_\sigma$ remains constant on each such set. 

We verify that, for each $t$ such that $\partial E_t$ is a $C^1$-hypersurface, each connected component $\Pi_t$  of $\partial E_t$ intersects $\partial \Omega$.  Indeed, if not, then the Alexander duality theorem (\cite{massey}) implies that $\Pi_t \cup \partial \Omega$ admits a decomposition of the form $\Pi_t \cup \partial \Omega = O_1 \cup O_2\cup(\R^n/\overline{\Omega})$ with $O_i$ open and connected.   We claim that $(\partial O_1 \cup \partial O_2)\cap \Pi_t \neq \emptyset$.  Indeed, were this not the case then $\partial O_i \cap \partial \Omega \neq \emptyset$ for $i=1,2$ and we could again apply the Alexander duality theorem to obtain the decomposition of $\R^n/\Pi_t$ into a union of two bounded and unbounded open, connected regions, $\Pi_t^c=O_b\cup O_u$.  Since, by the contrary assumption  $O_1 \cup O_2 \cup \Omega^c$ is connected and unbounded we have that $O_1 \cup O_2 \cup \Omega^c \subset O_u$ and thus, on taking complements, $O_b \subset (O_1 \cup O_2 \cup \Omega^c)^c \subset \Pi_t$. This contradiction shows that $(\partial O_1 \cup \partial O_2)\cap \Pi_t \neq \emptyset$, as claimed.  It thereby follows that $\Pi_t \cap \partial \Omega \neq \emptyset$. In other words, connected components of almost all the level sets $\partial E_t$ reach the boundary $\partial \Omega$. 

Now define $G\coloneqq \{t\in \R_+: \tilde{u}\left.\right|_{\partial E_t}={u_\sigma}\left.\right|_{\partial E_t}\}\subset \R_+$ and suppose there exists a ball $B \subset \Omega$ whose closure is contained in $\Omega$ and such that $\overline{B} \cap \{\tilde{u} \in G \}=\emptyset$.  If so then since $|\nabla \tilde{u}|\left.\right|_{B}\neq 0$  $\mu - a.e.$,  $\tilde{u}$ must map $\overline{B}$ to the closed interval $[a,b]$.  However $\mu([a,b]) >0$ whereas $\mu(G^c\cap \Omega)=0$, a contradiction.  We conclude from this that the union of all preimages of $G$ under $u$ and $\tilde{u}$ is dense in $\overline{\Omega}$ and therefore  $u=\tilde{u}$ $\mu-a.e.$, as was to be demonstrated. 

Finally, with $J=-c\sigma_0 \nabla u_\sigma$ we have $(\sigma_0^{-1} J \cdot J)^{\frac{1}{2}}=(c^2 \sigma_0 \nabla u_\sigma \cdot \nabla u_\sigma)^{\frac{1}{2}}$ $\mu - a.e.$. This gives the desired formula for $c(\x)$.
}

\end{proof}
\end{theorem}

\subsection{Equipotential Sets are Minimal Surfaces in a Riemannian metric Determined from the Data}

We close this section with some interesting geometrical results about the level sets of solutions to (\ref{bvp1}). Given $\sigma_0$ and the magnitude $|\J|_{\sigma_0^{-1}}$ of the current, we define a Riemannian metric on $\Omega$ and show that the level sets of the corresponding potential function have zero mean curvature in this metric.  In section \ref{geometric} we will prove the stronger statement that these equipotential sets are in fact area minimizing. These are generalizations to anisotropic conductivites of results proved in \cite{alex2, alex4} for the isotropic case.

As is customary, we denote $|A| \coloneqq \det A$ for $A \in Mat(\R,n)$ (which should not be mistaken  for the norm $|V|_{\sigma_0}$ of a vector field $V$, as we hope shall be clear from the context).
\begin{proposition}
\label{mean}
{Let $\Omega \subset \R^n$, $n \geq 2$ be a domain with Lipschitz boundary and $u \in C^{1, \alpha}(\bar{\Omega})$, $\alpha >0$.  Assume the conductivity $\bsigma$ is of the form $c(\x) \bsigma_0(\x)$ for $c,  \bsigma_0 \in C^\alpha(\Omega)$ with $\bsigma_0$ a known positive-definite matrix-valued function and that $|\nabla u|, c(\x) >0$ $\mu$-a.e. where $u$ is the potential corresponding to the conductivity $\sigma$ and current density $J$ via $J=-\sigma \nabla u$. 

Define the following Riemannian metric $g_{ij}$ on $\Omega$: \be
g_{ij}\coloneqq (|\bsigma_0||\J|_{\sigma_0^{-1}}^{2})^{\frac{1}{n-1}} (\sigma_0^{-1})_{ij}.
\ee
Then inside $\Omega$ one has that
\[\nabla \cdot (\sqrt{|g|}\frac{g^{ij}\nabla_i u}{||g^{-1}\nabla u||_g})=0.
\]
}
\begin{proof}
We begin by noticing that $|\bsigma_0|^{\frac{1}{n-1}}|\J|_{\sigma_0^{-1}}^{\frac{2}{n-1}} \bsigma_0^{-1}=c^{1+\frac{1}{n-1}-\frac{n}{n-1}}\{|\bsigma|(\bsigma \nabla u \cdot \nabla u)\}^{\frac{1}{n-1}}\bsigma^{-1}$ whereby, with the above choice of $g_{ij}$ one has that 
\[g^{-1}=\{|\bsigma|(\bsigma \nabla u \cdot \nabla u)\}^{\frac{1}{1-n}}\bsigma \]
Defining $m(x)\coloneqq |\bsigma|(\bsigma \nabla u \cdot \nabla u)$ it immediately follows that $|g|=\frac{m^{\frac{n}{n-1}}}{|\bsigma|}$.  Since $||g^{-1} \nabla u||^2_g = \{(g^{-1} \nabla u) \cdot g (g^{-1} \nabla u)\}^2$ we have $||g^{-1} \nabla u||_g= \sqrt{(g^{-1}\nabla u) \cdot \nabla u }$.  Then 
\begin{align}
\nabla \cdot (\sqrt{|g|}\frac{g^{ij}\nabla_i u}{||g^{-1}\nabla u||_g})&=\nabla \cdot(\frac{m^{\frac{n+1}{2(n-1)}-\frac{1}{n-1}} \bsigma \nabla u}{\sqrt{|\bsigma| \bsigma \nabla u \cdot \nabla u }} )\nonumber\\
&=\nabla \cdot(\frac{\sqrt{m(x)}\bsigma \nabla u}{\sqrt{m(x)}})\nonumber
\end{align}
It  follows from the fact that $u$ solves the conductivity equation that 
\[\nabla \cdot (\sqrt{|g|}\frac{g^{ij}\nabla_i u}{||g^{-1}\nabla u||_g})=0.\]

\end{proof}
\end{proposition}
The above result immediately implies the following.
\begin{corollary}{Suppose that $u, c, \sigma_0$ are as is in proposition (\ref{mean}).  Then the level sets of $u$, $u^{-1}(\lambda) \coloneqq \{u=\lambda\} \cap \overline{\Omega}$  are surfaces of zero mean curvature in the metric 
\[ g_{ij} = ( |\bsigma_0| |\J|_{\sigma_0^{-1}}^{2})^{\frac{1}{n-1}} (\sigma_0^{-1})_{ij}.\]
%Moreover, level sets of $u$, are critical hypersurfaces of the area functional
%\[A(\Sigma)\coloneqq \int_\Sigma |J|_{\sigma_0^{-1}} d\mu_{n-1}, \qquad \Sigma \subset \Omega \]
%where $d\mu_{n-1}$ is the surface measure on $n-1$ dimensional subspaces $\Sigma \subset \Omega$ and where $(g)^{\frac{1}{2}}$ is a choice of square root of $g_{ij}$.
}

\begin{proof} As in the preceding proof of Theorem \ref{mainthm} the level sets $u^{-1}(\lambda)$ are $C^1$-hypersurfaces for $\mu-$a.e. $\lambda$.  The unit vector $\eta \coloneqq \frac{g^{-1}\nabla u}{|\nabla u|_g}$ is $g$-orthogonal to each such level set $u^{-1}(\lambda)$ since if $\xi \in T_x u^{-1}(\lambda)$ then $\xi \cdot \nabla u\left. \right|_{x \in u^{-1}(\lambda)}=0$ and therefore 
\begin{align}
g(\xi, \eta)&=\frac{g(\xi, g^{-1}\nabla u)}{|\nabla u|_g}\nonumber \\
&= \frac{\xi \cdot \nabla u}{|\nabla u|_g}\nonumber \\
&= 0\nonumber 
\end{align}
Thus $\pm \eta$ are unit normals to the hypersurfaces $u^{-1}(\lambda)$ in the $g_{ij}$ metric.  Since the mean curvature of a surface with unit normal $n$ is $H=div_g(n)$, with $div_g$ being the metric divergence, we conclude from Proposition (\ref{mean}) that when $u$ satisfies the conductivity equation, we have $H=0$.
  
%
%Next, from Theorem $2.1$ in \cite{alex4} the level sets, $\tilde{u}_\lambda$ of $C^1(\Omega)$-solutions to the isotropic conductivity equation
%\[\nabla \cdot (\frac{|\tilde{J}|\nabla \tilde{u}}{|\nabla \tilde{u}|})=0, \qquad \tilde{J}=-\tilde{\sigma}\nabla \tilde{u}\neq 0   \]
%are critical hypersurfaces of $\int_{\Sigma}\sqrt{\tilde{J}\cdot \tilde{J}} dS$. Define $\tilde{J}\coloneqq \gamma J$ where $\gamma$ is a choice of $\sqrt{\sigma^{-1}_0}$.  Then $|J|_{\sigma_0^{-1}}=\sqrt{\tilde{J}\cdot \tilde{J}}$ and the critical surfaces corresponding to level sets of $\tilde{u}$, $\Sigma_c$, of $A(\cdot )$ have normal vectors that are orthogonal, with respect to the the Euclidean metric, to $\gamma^{-1}\nabla u$, since $\nabla \tilde{u}=\frac{c}{\tilde{\sigma}}\gamma^{-1}\nabla u$.  But, observe that $\sqrt{g^{-1}}$ and $\gamma^{-1}$ are conformally equivalent.  Therefore, if $\xi \in T_x\Sigma_c$ is a vector tangent to such a critical hypersurface  we have
%\begin{align}
%0&=\xi \cdot \nabla \tilde{u}\left.\right|_{\Sigma_c}\nonumber \\
%&= g(\sqrt{g^{-1}}\xi, \sqrt{g^{-1}}\nabla \tilde{u})\left.\right|_{\Sigma_c} \nonumber \\
%\end{align}
%so that 
%\[(\sqrt{g}\left.\right|_{\Sigma_c} \xi)\cdot \nabla u=0\]
%
\end{proof}
\end{corollary}

\section{Preliminaries for the General Case} 
\label{prelim}
In this section we prepare to expand upon the results in the preceding section by considering the conductivity equation over domains which may contain insulating or perfectly conducting inclusions, i.e. regions of infinite or zero conductivity, respectively.   We shall give the appropriate reformulation of  the forward problem \eqref{bvp1} in this setting.  We also discuss integration by parts and coarea formulae   for spaces of bounded weighted variation which will play a key role in our main general uniqueness result.

\subsection{Weighted Total Variation}
We start by presenting some needed preliminary results about functions of bounded weighted variation.  We will always use the notation $\chi_A(x)$ to denote the characteristic function of a set $A$.   Often, we will abbreviate vectors and matrices in component form, in addition, as earlier, we will employ the Einstein summation convention of implied summation over repeated upper and lower indices wherever appropriate.  
 
Let $\Omega \subset \R^n$ be a bounded open set
with connected Lipschitz boundary and let $a\in L^{\infty}
(\Omega)$ be non-negative.  While the function $a$ is now allowed to vanish, we require that its zero set  $S\coloneqq\{x\in \bar{\Omega}: a(x)=0\}$ always be assumed to satisfy the following structural hypothesis
\be
\label{hyp}
S\coloneqq O\cup \Gamma,
\ee
where $\Gamma$ is a set of measure zero with at most countably many connected components, $\Haus^{n-1}(\partial \Omega \cap S)=0$, and where $O$ is a mutually disjoint union of finitely many $C^1$-diffeomorphic  images of the unit ball, possibly empty. 

We recall the space of functions of locally bounded variation on a subset $\Sigma \subset \R^n$ is given by \[BV_{loc}(\Sigma)\coloneqq \{u\in Lip(\Sigma), \int_K |Du|<\infty, \quad \forall K\subset \Sigma,\quad K \text{ compact} \} \]
where $Du$ is the distributional gradient of $u$.  This generalizes the space $BV(\Omega)$, the space of all $L^1(\Omega)$ functions with bounded \textit{total variation} of the distributional gradient, i.e. those functions satisfying
\[\int_\Omega |Du|<\infty \]

Let $\sigma_0 \in C^0 (\Omega , Mat(\R,n))$ be a symmetric positive definite matrix with components $(\sigma_0)_{ij}$ satisfying 
\[m |\xi|^2 \leq \sum_{i,j=1}^n(\sigma_0)_{ij}(x) \xi^i \xi^j \leq M |\xi|^2 \ \ \ \ \forall x\in \Omega \setminus S, \ \ \forall \xi \in \R^n,\]
for constants $0<m, M< \infty$.  We then denote by $\phi(x,\xi)$ the following function
\begin{equation}
\varphi(x,\xi)=a(x)(\sum_{i,j=1}^n(\sigma_0)_{ij}\xi^i \xi^j  )^{\frac{1}{2}}. 
\end{equation}
For $u \in BV_{loc}(\Omega \setminus S)$ we define the \textit{weighted total variation} of $u$, with respect to $\varphi$, in $\Omega$ as 
\begin{equation}
\label{BV-def}
\int_{\Omega} |Du|_{\varphi}=\sup_{{B} \in \mathfrak{B}_a} \int_{\Omega} u  \nabla \cdot {B}  \ \ d\mu,
\end{equation}
where 
\begin{equation*}
\mathfrak{B}_a=\{ B \in L_c^{\infty}(\Omega,\R^n) :\quad   \nabla \cdot {B} \in L^{n}(\Omega) \quad \hbox{and} \quad |B|_{\sigma_0^{-1}}\leq a(x) \quad  \hbox{a.e. in} \quad \Omega\}.
\end{equation*}
and $L_c^{\infty}(\Omega,\R^n)$ is the space of vector fields of compact support in $\Omega$ whose components are in $L^\infty(\Omega)$.

It is a straightforward consequence of the definition (\ref{BV-def}) that $\int_{\Omega} |Du|_{\varphi}$ is $L_{loc}^{\frac{n}{n-1}}(\Omega)-$lower semi-continuous.  
It was shown in \cite{AB} by Amar and Bellettini that for any $u\in BV(\Omega)$,  one has the following integral representation formula for the weighted total variation appearing in equation (\ref{BV-def}),
\begin{equation}\label{BV-defAB}
\int_{\Omega}|Du|_{\varphi}=\int_{\Omega} h(x,v^u)|Du|
\end{equation}
where, in the above, 
\begin{equation}
h(x,v^u)\coloneqq (|Du|-\hbox{ess} \sup_{\mathbf{B}\in \mathfrak{B}} ({B} \cdot v^u))(x) \qquad \qquad |Du|-a.e. \ \ x\in \Omega,
\end{equation}
and $v^u$ denotes the vectorial Radon-Nikodym derivative $v^{u}(x)=\frac{d\, Du}{d\, |Du|}$.
Note that the right-hand side of equation \eqref{BV-defAB} makes sense, as $v^u$ is $|Du|$-measurable,
and hence $h(x, v^u(x))$ is as well.  In particular, it can be shown (viz. \cite{AB} Prop. 7.1) that if $a$ and $\sigma_0$ are continuous in $\Omega$, then one has 

\begin{equation}
\label{useful}
h(x,v^{u})=a(x)\left( \sum_{i,j=1}^{n}\sigma_0^{ij} v_i^{u}v_j^{u} \right)^{1/2},  \qquad |Du|-a.e. \text{ in } \Omega 
\end{equation}
for every Borel set $\Omega$ and $u\in BV(\Omega)$.

Following \cite{Al} and \cite{An}, we let
\[ X\coloneqq \{B\in L^{\infty}(\Omega,\R^n): \text{div } {B} \in L^n(\Omega)\}. \]
As proven in \cite{An}, Theorem 1.2, if $\nu_\Omega$ denotes the outer unit normal vector to $\partial \Omega$, then for every $B \in X$ there exists a unique function $[B \cdot \nu_{\Omega}] \in L^{\infty}_{\mathcal{H}^{n-1}}(\partial \Omega)$ such that   
\begin{equation}
\int_{\partial \Omega} [B \cdot v_{\Omega}]u d\mathcal{H}^{n-1}=\int_{\Omega} u \nabla \cdot B d\mu +\int_{\Omega}  B \cdot \nabla u d\mu, \ \quad \forall u \in C^1(\bar{\Omega}). 
\end{equation}
Moreover, for $u\in BV(\Omega)$ each such $B\in L^{\infty}(\Omega, \R^n)$ with $\nabla \cdot {B} \in L^n(\Omega)$ gives rise to a Radon measure on $\Omega$, denoted $({B} \cdot Du)$, satisfying the following
\begin{equation}\label{IBP0}
\int_{\partial \Omega} [{B}\cdot v_{\Omega}]u d\mathcal{H}^{n-1}=\int_{\Omega} u \nabla \cdot {B} d\mu +\int_{\Omega}  ({B} \cdot D u), \qquad \forall u \in BV(\Omega),
\end{equation}
We refer the interested reader to \cite{Al, An} for a proof. 

We shall need the following lemma, a proof of which follows from (\ref{IBP0}), and the fact that
\[BV_{loc}(\Omega \setminus \bar{S}) \cap L^{\infty}(\Omega) \subset BV(\Omega)\]
which can be easily verified. 

\begin{lemma}\label{IBP} Let $S$ be as defined in (\ref{hyp}). Then 
\begin{equation}
\int_{\partial \Omega} [{B}\cdot v_{\Omega}]u d\mathcal{H}^{n-1}=\int_{\Omega} u \nabla \cdot {B} d\mu +\int_{\Omega}  ({B} \cdot D u) 
\end{equation}
for all $u \in BV_{loc}(\Omega \setminus \overline{S}) \cap L^{\infty}(\Omega)$. 

\end{lemma}

We conclude with a useful co-area formula for functions of bounded weighted total variation. Details can be found in \cite{AB}.
\begin{theorem}[Generalized Co-Area Formula]{Let $u \in BV(\Omega)$ and suppose $\Haus^{n-1}(\Omega \cap \{u=s\})<\infty$ holds for all $s \in \R$. Let $P_\phi(A,\Omega)$ denote the perimeter of the set $A \subset \Omega$ given by 
\[P_\phi(A,\Omega) \coloneqq \int_\Omega |D\chi_A|_\phi \]  
  Then
\[\int_\Omega |Du|_\phi =\int_\R P_\phi(\{u>s\}, \Omega)ds \]
}
\end{theorem}
We note that this may, on using the representation formula \eqref{BV-defAB},  be recast as
\be
\label{coarea}
\int_\Omega|Du|_\phi = \int_\R \int_{\Omega \cap \partial^* \{u(x)>s \}}h(x,v^s)d\Haus^{n-1}(x)ds
\ee
where $v^s$ is a unit outer-oriented normal vector to $\Omega \cap \partial^*\{u(x)>s \}$.

\subsection{Modeling Regions with Zero or Infinite Conductivity}

Here we discuss how to formulate a suitable version of the conductivity equation \eqref{bvp1} in the presence of inclusions of infinite and/or zero conductivity. 

Let $O_{\infty}$ be an open subset of $\Omega$ satisfying
$\overline{O}_{\infty}\subset\Omega$, meant to model \emph{perfectly conducting
inclusions}, and $O_{0}$ be an open subset of $\Omega$ with
$\overline{O}_0\subset\Omega$, meant to model \emph{insulating
inclusions}.  We assume
$\overline{O}_\infty\cap\overline{O}_0=\emptyset$,  $\Omega\setminus
\overline{O_\infty\cup O_0}$ is connected, and the boundaries $\partial O_\infty$,
$\partial O_0$ are piecewise $C^{1,\alpha}$ for $\alpha >0$. We also assume that $O_0$ is a mutually disjoint union of finitely many $C^1-$ diffeomorphic images of the unit ball, possibly empty. In addition, in two dimensions we require that $O_0$ has at most one such component.

Let $\sigma^{jk}$ and $\tilde{\sigma}^{jk}$ be symmetric positive definite matrix functions. For  $k>0$ consider
the conductivity problem
\begin{equation}\left\{ \begin{array}{ll}
\partial_{x_j}\left\lbrace \left[ (k \tilde{\sigma}^{ij}-\sigma^{ij})\chi_{O_\infty}+\sigma^{ij}\right] \partial_{x_i}u_k \right\rbrace=0, \ \ \hbox{in} \ \  \Omega \setminus \overline{O}_0\\
\frac{\partial u_k}{\partial \nu}=0 \ \  \hbox{on}\ \ \partial O_0, \\
 u_{k}|_{\partial\Omega}=f.
\end{array} \right.
\end{equation}

The perfectly conducting inclusions occur in the limiting case $k\to
\infty$. The limiting solution is the unique solution to the
problem:
\begin{equation}\left\{ \begin{array}{ll}
\partial_{x_j}\left(\sigma^{ij}\partial_{x_i}u \right)=0, \ \ &\hbox{in} \ \  \Omega \setminus \overline{O_0 \cup O_\infty}\\
%\nabla\cdot \sigma \nabla u_0=0,&\mbox{in}\,V\\
\nabla u=0, &\mbox{in} \ \ O_\infty\\
u|_+=u|_-,&\mbox{on}\ \ \partial (O_0\cup O_\infty)\\
\int_{\partial O_{\infty}^n}\mathbf{\sigma}\frac{\partial u}{\partial \nu_n}|_{+}dS=0,& n=1,2,...\\
\frac{\partial u}{\partial \nu}|_{+}=0,&\mbox{on}\;\partial O_0\\
u|_{\partial\Omega}=f,
\end{array}\label{pde_inclusions} \right.
\end{equation}
 (see the Appendix for more details), where  $O_\infty=\cup_{n=1}^\infty O_\infty^{n}$ is a  partition of $O_\infty$ into connected
components. Here , as in the rest of the paper, $\nu$ is the outward unit normal vector and the subscripts $\pm$ indicate the limit taken from the outside and inside the domain, respectively.

\begin{remark}For Lipschitz continuous conductivities in any dimension $n\geq 2$,
or for essentially bounded conductivities in two dimensions, the
solutions of the conductivity equation satisfy the unique
continuation property (see, \cite{bersJohnSchechter} and references
therein). Consequently the insulated (and possibly perfectly
conducting) inclusions are the only open sets on which the interior
data $|J|_{\sigma_0^{-1}}$ may vanish identically. However, in three dimensions or
higher it is possible to have a H\"{o}lder continuous $\sigma$ and
boundary data $f$ that yield $u\equiv constant$ in a proper open
subset $O_s \subsetneq \Omega$, see \cite{plis,martio}. We call such
regions $O_s$ {\em singular inclusions}. On the other hand, we will not use Ohm's law
in the classical sense inside perfect conductors: the current $J$ inside
perfectly conducting inclusions is not necessarily zero whereas
$\nabla u\equiv 0$ within such regions (see \cite{Aetall, Ip}).
\end{remark}

\section{Uniqueness and Determining the Conformal Class}
\label{inclusions}
From now on we assume that $\sigma \in C^{\alpha}(\Omega , Mat(\R,n))$ for $\alpha >0$ and satisfies
\begin{equation}\label{sig2}
\sigma(\x)=c(\x)\sigma_0(\x),
\end{equation}
where $c(\x) \in L^\infty_+(\Omega \setminus (O_\infty\cup O_0))$ is a real, scalar-valued function, bounded away from zero and finite on $\Omega \setminus (O_\infty\cup O_0)$ and where $\sigma_0 \in  C^{0}(\Omega, Mat(\R,n))$ is symmetric, positive definite.  

We will prove that the shape and locations of the perfectly conducting and insulating inclusions and the conductivity $\sigma$ outside of the inclusions are determined from knowledge of the boundary voltage $f$, $\sigma_0$ and 
\[a=\sqrt{\sigma^{-1}_0 {J} \cdot {J}}=|J|_{\sigma_0^{-1}} \ \ \hbox{in}\ \ \Omega \] 
where $J$ is the current density vector field generated by imposing the voltage $f$ at $\partial \Omega$. To formulate our results, we first need to extend the notion of \emph{admissibility}. 

\begin{definition} \label{def} {A triplet of functions $(f, \sigma_0, a)\in H^{1/2}(\partial \Omega)\times C^{0}(\Omega, Mat(\R,n)) \times
L^{\infty}(\Omega)$ is called \emph{admissible} if there exists a matrix valued function $\sigma$ satisfying ($\ref{sig2}$)  and a divergence free vector field $\J \in (L^{\infty}(\Omega))^n$ such that the following three statements hold.\\
(i) \begin{equation}
a= |J|_{\sigma_0^{-1}}  \ \ \hbox{in} \ \ \Omega.
\end{equation}
(ii) The vector field  $\J$ satisfies 
\begin{eqnarray*}\label{rel}
J=\left\{ \begin{array}{ll}
-\sigma \nabla u \ \  &\hbox{in}\ \ \Omega \setminus (O_\infty \cup O_0).\\
0 \ \ &\hbox{in} \ \ O_0,
\end{array} \right.
\end{eqnarray*}
where  $u$ is the corresponding solution of  ($\ref{pde_inclusions}$). \\
(iii) The set of zeros of $a$ can be decomposed as follows
\begin{equation}\label{zeroes }
\{x\in \Omega:\, a(x)=0\}\cap (\Omega\backslash \overline{O}_\infty)= O_0\cup
\overline{O}_s \cup \Gamma,
\end{equation}
where $O_s$ is an open set (possibly empty), $\Gamma$ is a
Lebesgue-negligible set, and $\overline{\Gamma}$ has empty interior. }
\end{definition}
We are now ready to state our main uniqueness results.

\begin{theorem} \label{unique}
Let $\Omega \subset \R^n$, $n\geq 2$, be a domain with connected
Lipschitz boundary and let $(f, \sigma_0, a)\in H^{1/2}(\partial \Omega)\times C^{0}(\Omega , Mat(\R,n)) \times
L^{\infty}(\Omega)$ be an admissible triplet generated by an
unknown  conductivity $\sigma$.  Define $\phi(x,\xi)= a(x)|\xi|_{\sigma_0}$ on $\Omega \setminus \bar{S}$. 

Then the potential $u$, solving \eqref{pde_inclusions}, is a minimizer of the problem
\begin{equation}\label{min_prob}
u=\hbox{argmin} \{\int_{\Omega} |Dv|_\phi: v \in BV(\Omega \setminus \bar{S}) \ \ \hbox{and} 
\ \ v|_{\partial \Omega}=f\},
\end{equation}
and, if $\bar{u}$ is another minimizer of the above problem, then
$\bar{u}=u$ in $\Omega \backslash \{|\J|=0\}$.

Moreover the zero-sets of $a$ and $|\nabla u |$ can be decomposed
as follows
\[\{x\in \Omega,a(x)=0\}\cup \{x\in \Omega:  \nabla u=0\}=:Z\cup \Gamma,\]
where $\Gamma$ has measure zero and $Z=O_\infty\cup O_0\cup O_s$ is open.  

Consequently \[\sigma=\frac{a}{|\nabla u|_{\sigma_0}} \sigma_0 \in
C^{\alpha}(\Omega\setminus \overline {Z})\] is the unique conductivity generating the admissible data triplet $(f, \sigma_0, a)$. 
\end{theorem}

\begin{remark} \label{rem1}The above theorem allows us to identify the potential $u$
and the conductivity $\sigma$ outside the open set $Z=O_\infty\cup O_0 \cup
S$. To determine if an open connected
component $O$ of $Z$ is a perfectly conducting inclusion, an
insulating inclusion, or a singular inclusion we proceed as follows:
\begin{itemize}
\item If $\nabla u\equiv 0$ in $O$ and $a(x)\neq 0$ for
some $x \in O$, then $O$ is a perfectly conducting inclusion.
\item If $a \equiv 0$ in $O$ and $ u\not\equiv constant$ on $\partial O$,
then $O$ is an insulating inclusion.
\item If $a \equiv0$ in $O$, $u = constant$ on $\partial O$, and $a$ is not
$C^{\alpha}$ at $x$ for some $x \in O$, then $O$ is either an
insulating inclusion or a perfectly conducting inclusion.
\item If $a\equiv0$, $u=constant$ on $\partial O$, and $a\in C^{\alpha}(\partial
O)$, then the knowledge of the magnitude of the current $(f, \sigma_0, a)$  is not
enough to determine the type of the inclusion $O$.
\end{itemize}
\end{remark}
\begin{proof}
%[{\bf Proof of Theorem \ref{unique}:}]  
Suppose that $\bar{u}\in BV_{loc}(\Omega \setminus \bar{S})$. First note that for every $x\in \Omega\setminus (\overline{O_\infty \cup O_0})$ there exists $\epsilon>0$ such that $B(x, 2\epsilon) \subset \Omega$ and 
\[\int_{B(x,  \epsilon)} h(x,v^{\bar{u}})|D \bar{u}|\geq -\int_{B(x,  \epsilon)} J \cdot v^{\ubar}|D\ubar|,\]
where $J$ is the current density vector field described in definition (\ref{def}). Therefore
\[h(x,v^{\ubar})\geq  - J\cdot v^{\ubar}, \ \ \ \ |D\ubar|-a.e. \ \ \hbox{in} \ \ \Omega\setminus (\overline{O_\infty \cup O_0}).\]
Thus, on using Lemma \ref{IBP} and the fact that the current density is divergence-free away from the inclusions we have
\begin{align}
\int_{\Omega \setminus \overline{(O_0 \cup O_{\infty})}} |D\bar{u}|_{\varphi}&=\int_{\Omega \setminus \overline{(O_0 \cup O_{\infty})}}h(x, v^{\bar{u}})|D\ubar| \nonumber  \\
& \geq - \int_{\Omega \setminus \overline{(O_0 \cup O_{\infty})}} J \cdot v^{\bar{u}} |D\ubar|\nonumber \\
&=- \int_{\Omega \setminus \overline{(O_0 \cup O_{\infty})}} J\cdot D\ubar \nonumber \\
&=-  \int_{\partial \Omega \setminus \overline{(O_0 \cup O_{\infty})}} J\cdot \nu f d\mathcal{H}^{n-1}\nonumber \\
&= \int_{\Omega \setminus \overline{(O_0 \cup O_{\infty})}} |Du|_{\varphi}. \nonumber 
\end{align}
 Hence 
\begin{equation}\label{h}
h(x,v^{\bar{u}})=-J\cdot v^{\bar{u}},\ \ |D \bar{u}|-a.e. \ \ \hbox{in} \ \ \Omega \setminus \overline{(O_0 \cup O_{\infty})}
\end{equation}
Then again, since $a$ is continuous in $\Omega \setminus (\overline{O_0 \cup O_{\infty}})$ equation \eqref{useful} gives  
\[ h(x,v^{\bar{u}})=a(x)\left( \sum_{i,j=1}^{n}\sigma_0^{ij} v_i^{\bar{u}}v_j^{\bar{u}} \right)^{1/2} \ \  |D \bar{u}|-a.e. \ \ \hbox{in} \ \  \Omega \setminus (\overline{O_0 \cup O_{\infty}}).\]  
But then, on $\Omega \setminus (\overline{O_0 \cup O_{\infty}})$, we have 
\begin{eqnarray*}
 h(x,v^{\bar{u}})&=& a(x)\left( \sum_{i,j=1}^{n}\sigma_0^{ij} v_i^{\bar{u}}v_j^{\bar{u}} \right)^{1/2}\\ 
 &=& c(x)|\nabla u|_{\sigma_0}| v^{\bar{u}}|_{\sigma_0}\\
 &\geq & c(x)|(\nabla u, v^{ \bar{u}})_{\sigma_0} |\\
 &\geq & \sigma \nabla u \cdot v^{\bar{u}}\\
 &=&-J \cdot v^{\bar{u}}. 
\end{eqnarray*}
Thus it follows from (\ref{h}) that
\[\frac{J}{|J|} =\frac{\nabla \ubar}{|\nabla \ubar |}=v^{\bar{u}}, \ \ |D \bar{u}|-a.e. \ \ \hbox{in} \ \ \Omega \setminus \overline{(O_0 \cup O_{\infty})}.\]
It follows from an argument similar to that of Theorem 3.5 in \cite{MNT12} that $u=\bar{u}$ a.e. in $\Omega$.  
\end{proof}

\section{Geometrical Properties of Equipotential Sets}
\label{geometric}
In this section we prove the area-minimizing property of the equipotential sets $u^{-1}(\lambda) \coloneqq \overline{\Omega} \cap \{x; u(x)=\lambda\}$ for solutions $u(x)$ of the equation  (\ref{pde_inclusions}). This generalizes results in \cite{alex2} and \cite{alex5}. The basic argument goes back to (\cite{bombieri}). We shall need the following lemmas, which require the generalized coarea formula (\ref{coarea}). 

Define the functional $\mathcal{F}$ by 
\[ \mathcal{F}[u]\coloneqq \int_{\Omega}|D v|_\phi \]
for the positive definite, matrix-valued $\sigma_0 \in C^{1,\alpha}(\overline{\Omega},Mat(\R,n))$. As well, for a given function $w$ on $\Omega$, we introduce the notation 
\[w_{\lambda, \eps}\coloneqq \frac{1}{\eps}\min\{\eps,\max{w-\lambda, 0}\} \]
for fixed $\lambda$ and $\eps>0$.  For functions $w$ on $\partial \Omega$, $w_{\lambda, \eps}$ denotes the trace at the boundary of such functions.

\begin{lemma}
{Let  $a \in Lip(\Omega)$, fix $\lambda \in \R$ and let $\eps>0$.  Then if 
\[u \in \underset{v \left. \right|_{\partial \Omega}=f}{argmin}\{\int_\Omega |Dv|_\phi; \ v \in W^{1,1}_+(\Omega)\} \] 
Then also
\[u_{\lambda, \eps} \in \underset{v \left. \right|_{\partial \Omega}=f_{\lambda, \eps}}{argmin}\{\int_\Omega |Dv|_{\phi}; \ v \in W^{1,1}_+(\Omega)\} \] 
where $ u_{\lambda, \eps}=\frac{1}{\eps}\min\{\eps,\max{u-\lambda, 0}\}$.
}
\begin{proof}
For $u \in W^{1,1}_+(\Omega)$ we have $v^u=\frac{\nabla u}{|\nabla u|}$ $\mu -a.e.$ Set $u_+=\max\{u-\lambda, 0\}$ and $u_-=u-u_+$. Then since $a$ and $\sigma_0$ are continuous we have, with $\phi(x,\xi)=a|\xi|_{\sigma_0}$, that 
\[ \int_{\Omega} |Du|_\phi=\int_{\Omega}a(x)(\sigma_0 v^u \cdot v^u)|\nabla u|=\mathcal{F}[u]\]
Then again, by (\ref{coarea}), we have 
\begin{align}
\mathcal{F}[u]&=\int_{t<\lambda }\int_{u^{-1}(t)}a(x)d\Haus^{n-1}(x)dt+\int_{t>\lambda}\int_{u^{-1}(t)}a(x)d\Haus^{n-1}(x)dt\nonumber \\
&=\int_{t \in \R} \{\int_{u_+^{-1}(t)}a(x) \Haus^{n-1}(x)+ \int_{u_-^{-1}(t)}a(x)d\Haus^{n-1}(x)\}dt \nonumber \\
&=\mathcal{F}[u_+]+\mathcal{F}[u_-]
\end{align}
Suppose that $u \in \underset{v\left.\right|_{\partial \Omega}=f}{argmin} \{\mathcal{F}[v]; v \in W^{1,1}_+(\Omega)\}$, whereby one has with a test function $w\in Lip(\Omega)$ chosen such that $w\left.\right|_{\partial \Omega}=0$ that 
\begin{align}
\mathcal{F}[u_+]&=\mathcal{F}[u]-\mathcal{F}[u_-]\nonumber\\
&\leq \mathcal{F}[u+w]-\mathcal{F}[u_-]\nonumber \\
&\leq \mathcal{F}[u_+  + u_- - u_- + w] \nonumber
\end{align}
so that $u_+ \in \underset{v\left.\right|_{\partial \Omega}=f_+}{argmin} \{\mathcal{F}[v]; v \in W^{1,1}_+(\Omega)\}$ and, likewise $u_- \in \underset{v\left.\right|_{\partial \Omega}=f_-}{argmin} \{\mathcal{F}[v]; v \in W^{1,1}_+(\Omega)\}$ follows mutatis mutandis.  The lemma is immediate from two applications of this fact.
\end{proof}
\end{lemma}
\begin{lemma}{Let $a \in Lip(\Omega)$ and $u\in Lip(\Omega)$ be such that 
\[\{x;|\nabla u|=0\}=Z\cup L \]
where $Z$ is open, $a(\overline{Z})=0$, $\mu(L)=0$ and 
\[\frac{a \nabla u}{|\nabla u|} \in W^{1,1}(\Omega/ \overline{Z}) \]
Then 
\be
\lim_{\eps \downarrow 0}\mathcal{F}[u_{\lambda, \eps}]\label{limit}=\int_{u^{-1}(\lambda)}a(x)d\Haus^{n-1}(x)
\ee
holds for almost every $\lambda \in \R$.
}
\begin{proof}
The proof is similar to that of Lemma 4.3 in \cite{alex5}.
\end{proof}
\end{lemma}

We consider the functional
\be
\label{areafunc}
\mathcal{A}(\Sigma)=\int_{\Sigma}|J|_{\sigma_0^{-1}} dS
\ee
for $\Sigma \subset \Omega$ and $dS$ the induced Euclidean surface measure on the $n-1$ dimensional subset $\Sigma$.  We remark that when $\Omega$ is equipped with the data dependent metric $g_{ij} =(|\bsigma_0||\J|_{\sigma_0^{-1}}^{2})^{\frac{1}{n-1}} (\sigma_0^{-1})_{ij}$ discussed in Proposition \ref{mean} that the invariant volume form $\sqrt{|g|}dx^1 \wedge \cdots \wedge dx^n$ on $\Omega$ induces the invariant volume form $|J|_{\sigma_0^{-1}}dS(x)$ on smooth $n-1$ dimensional hypersurfaces $\Sigma$.  For this reason we refer to \eqref{areafunc} as an area functional.

We are ready to establish the main  result of this section, which says that equipotential hypersurfaces of solutions to ($\ref{pde_inclusions}$) are minimizers of  the area functional (\ref{areafunc}). 
\begin{theorem}
{Let $\Omega \subset \R^n$, $n\geq 2$, be a domain with connected
Lipschitz boundary and let $(f, \sigma_0, a)\in C^{2,\delta}(\partial \Omega)\times C^{1,\delta}(\Omega , Mat(\R,n)) \times
L^{2}(\overline{\Omega})$ be an admissible triplet generated by an
unknown  conductivity $\sigma \in C^{1,\delta}(\Omega \setminus (O_\infty\cup O_0)$ with $\delta \in (0,1)$.  Let $v\in C^{2,\delta}(\overline{\Omega})$ satisfy $v \left.\right|_{\partial \Omega}=f$ and
\[\{x;|\nabla v|=0 \}=Z_v \cap L_v, \quad a(\overline{Z_v})=0 \]
for $Z_v$ open and $\mu(L_v)=0$.  

Then, when $u$ is the solution to the BVP (\ref{pde_inclusions}) we have 
\be
\mathcal{A}(u^{-1}(\lambda))\leq \mathcal{A}(v^{-1}(\lambda))
\ee
for almost all $\lambda \in \R$. 
}
\begin{proof}
Notice that when $\lambda \notin Range (u)$, the result is immediate.  Otherwise since, under these hypotheses on the conductivity, $u$ satisfies the maximum principle and, by assumption $u$ and $v$ agree on $\partial \Omega$, $Range(u) \subset Range (v)$. If $\lambda \in Range(u)/\{u(\overline{Z})\cup v(\overline{Z_v})\}$ is arbitrary, $u(\overline{Z})\cup v(\overline{Z_v})$ is at most countable, and since $\nabla u$ and $\nabla v$ are both non-zero away from $\overline{Z}$and $\overline{Z_v}$ respectively, $u^{-1}(\lambda)$ and $v^{-1}(\lambda)$ are $C^2-$smooth oriented hypersurfaces whereon the Hausdorff measure coincides with the standard Lebesgue measure.  Thus $\Haus^{n}(u^{-1}(\lambda))=0$ for $a.e.$ $\lambda$ and, by assumption $\Haus^{n-1}(u^{-1}(\lambda) \cap \partial \Omega)=0$ for $a.e.$ $\lambda$ since $u$ extends continuously to the boundary, which has finite Hausdorff measure.

Since $u$ and $v$ agree on the boundary $u_{\lambda,\eps}\left.\right|_{\partial \Omega}=v_{\lambda,\eps}\left. \right|_{\partial \Omega}$ as well.  But then by Theorem \ref{unique} we have
\[ \mathcal{F}[u_{\lambda, \eps}] \leq \mathcal{F}[v_{\lambda, \eps}]\]
The result then follows from (\ref{limit}).% and the fact that $W^{1,1} \subset BV$ is a closed subset \cite{evans2}. 

\end{proof}
\end{theorem}

\section{Appendix: Perfectly conductive and insulating inclusions}
\label{appendix}
The results presented in this appendix give a precise definition, by a limiting procedure,  of  potentials corresponding to conductivities that can vanish or be infinite in certain regions. 
They are slight generalization of the ones in
\cite{baoLiYin} to include both perfectly conductive and insulating
inclusions.
%U=O_\infty and V=O_0

Let $O_\infty=\cup_{j=1}^\infty O^{j}_{\infty}$ be an open subset of $\Omega$ with
$\overline{O}_{\infty}\subset\Omega$ to model the union of the connected
components  $O^j_{\infty}$ ($j=1,2,...$) of \emph{perfectly conductive
inclusions}, and let $O_0$ be an open subset of $\Omega$ with
$\overline{O}_0\subset\Omega$ to model the union of all connected
\emph{insulating inclusions}. Let
 $\chi_{O_\infty}$ and $\chi_{O_0}$ be their corresponding characteristic
function. We assume that $\overline{O_\infty}\cap\overline{O_0}=\emptyset$,
 $\Omega\setminus \overline{O_{\infty}\cup O_0}$ is connected, and that the
boundaries $\partial O_\infty$, $\partial O_0$ are piecewise $C^{1,\alpha}$ for $\alpha >0$.
Let $\sigma_1\in C^{\alpha}(O_\infty, Mat(\R,n))$, and $\sigma\in
C^{\alpha}(\Omega\setminus \overline{O_{\infty}\cup O_0}, Mat(\R,n))$ be matrix-valued functions such that
\begin{equation*}
m |\xi|^2 \leq \sigma^{ij} \xi_i \xi_j \leq M |\xi|^2, \ \ m |\xi|^2 \leq \sigma_1^{ij} \xi_i \xi_j \leq M |\xi|^2
\end{equation*}
for  constants $0<m, M<\infty$ and all $\xi \in \R^n$. 

For each $0<k<1$ consider the conductivity problem
\begin{equation}
\nabla\cdot(\chi_U(\frac{1}{k} \sigma_1-\sigma)+\sigma)\nabla u=0,
\qquad \frac{\partial u}{\partial \nu}=0 \ \hbox{on} \ \
\partial O_0, \ \ \hbox{and} \qquad u|_{\partial\Omega}=f.\label{condEQ}
\end{equation} 
The condition on $\partial O_{0}$ ensures that $O_0$ is insulating. 
It is well known that the problem (\ref{condEQ}) has a unique
solution $u_k\in H^1(\Omega)$ which also solves
\begin{equation}
\left\{ \begin{array}{ll}
\nabla\cdot \sigma \nabla u_k=0,&\mbox{in}\,\Omega\setminus\overline{O_{\infty}\cup O_0},\\
\nabla\cdot \sigma_1 \nabla u_k=0,&\mbox{in}\,O_{\infty},\\
%\nabla\cdot \sigma_2\nabla u_k=0,&\mbox{in}\,V,\\
u_k|_+=u_k|_-,&\mbox{on}\,\partial O_{\infty},\\
%u_k|_+=u_k|_-,&\mbox{on}\,\partial V,\\
\left. (\frac{1}{k} \sigma_1 \nabla u_k) \cdot \nu\right|_-=\left. (\sigma \nabla u_k) \cdot \nu \right|_+,&\mbox{on}\,\partial O_{\infty},\\
\left. \frac{\partial u_k}{\partial \nu}\right|_{+}=0,&\mbox{on}\,\partial O_0,\\
u_k|_{\partial\Omega}=f.
\end{array}\label{pde_k} \right.
\end{equation} Moreover, the energy functional
\begin{equation}
I_k[v]=\frac{1}{2k}\int_{O_\infty} |\nabla v|^2_{\sigma_1} dx+\frac{1}{2}\int_{\Omega\setminus \overline{U\cup
V}} |\nabla v|^2_ \sigma dx
\end{equation}
has a unique minimizer over the maps in $H^1(\Omega)$ with trace $f$
at $\partial\Omega$ which is the unique solution $u_k$ of
(\ref{pde_k}).

We shall show below that the limiting solution (with $k\to 0$) solves

\begin{equation}\left\{ \begin{array}{ll}
\nabla\cdot \sigma \nabla u_0=0,&\mbox{in}\,\Omega\setminus\overline{O_{\infty} \cup O_0},\\
%\nabla\cdot \sigma \nabla u_0=0,&\mbox{in}\,V\\
\nabla u_0=0, &\mbox{in} \ \ O_{\infty},\\
u_0|_+=u_0|_-,&\mbox{on}\ \ \partial O_{\infty},\\
%u_0|_+=u_0|_-, &\mbox{on} \, \partial U\\
%u_0|_+=u_0|_-, &\mbox{on} \, \partial V\\
\int_{\partial O_{\infty}^j} (\sigma \nabla u) \cdot \nu |_{+}ds=0,&j=1,2,...,\\
\frac{\partial u_0}{\partial \nu}|_{+}=0,&\mbox{on}\;\partial O_0,\\
u_0|_{\partial\Omega}=f,\\
%\nabla\cdot \sigma_2 \nabla u_0=0,&\mbox{in}\;V,\\
%u_0|_+=u_0|_-,&\mbox{on}\;\partial V.
\end{array}\label{pde_0} \right.
\end{equation}
By elliptic regularity $u_0 \in C^{1,\alpha}(\Omega \backslash O_\infty \cup
O_0)$ and for any $C^{1,\alpha}$ boundary portion $T$ of $\partial
(O_\infty \cup O_0)$, $u_0 \in C^{1,\alpha}((\Omega \backslash (O_\infty \cup O_0))\cup
T)$.

\begin{proposition}
The problem (\ref{pde_0}) has a unique solution in $H^1(\Omega)$.
This solution is the unique minimizer of the functional
\begin{equation}\label{I0}
I_0[v]=\frac{1}{2}\int_{\Omega\setminus \overline{O_{\infty} \cup
O_0}}|\nabla v|^2_{\sigma} dx,
\end{equation}over the set $A_0\coloneqq\{u\in H^1(\Omega\setminus\overline{O}_0);\,
u|_{\partial\Omega} =f,\, \nabla u=0\, \mbox{in}\, O_{\infty} \}$.
\end{proposition}
{\bf Proof:} Note that $A_0$ is weakly closed in
$H^1(\Omega\setminus \overline{O}_0)$. The functional $I_0$ is lower
semicontinuous and strictly convex and, as a consequence,  has a unique minimizer
$u_0^*$ in $A_0$.

First we show that $u_0^*$ is a solution of the BVP (\ref{pde_0}). Since
$u_0^*$ minimizes (\ref{I0}), we have
\begin{equation}\label{weaksolution}
0=\int_{\Omega\setminus \overline{O_0\cup O_{\infty}}}\sigma \nabla
u^*_0\cdot\nabla\phi dx,
\end{equation}
for all $\phi\in H^1(\Omega\setminus \overline{O_0})$, with
$\phi|_{\partial\Omega}=0$, and $\nabla \phi =0 $ in $O_{\infty}$. In
particular, if $\phi\in H^1_0(\Omega\setminus \overline{O}_0)$, we get
$\int_{\Omega\setminus \overline{O_{\infty}\cup O_0}}(\nabla \cdot\sigma\nabla
u_0^*)\phi dx=0$ and thus $u_0^*$ solves the conductivity equation
in (\ref{pde_0}). If we choose $\phi\in H^1(\Omega\setminus
\overline{O}_0)$, with $\phi|_{\partial\Omega}=0$, and $\phi \equiv 0 $ in
$O_{\infty}$, from Green's formula applied to (\ref{weaksolution}), we get
$\int_{\partial O_0} (\sigma \nabla u_0^*) \cdot \nu |_+\phi=0, \,\forall \,\phi|_{\partial
O_0}\in H^{1/2}(\partial O_0),$ or equivalently, $ \frac{\partial u^*_0}{\partial \nu} |_{\partial O_0}=0 $. If we choose $\phi_j\in H^1_0(\Omega\setminus \overline{O}_0)$ with
$\phi_j\equiv 1$ in the connected component $O^j_{\infty}$ of $O_{\infty}$ and
$\phi_j\equiv 0$ in $O_{\infty} \setminus O_{\infty}^j$, from Green's formula applied
to (\ref{weaksolution}) we obtain $\int_{\partial
O_{\infty}^j}(\sigma \nabla u_0^*) \cdot \nu |_+=0$. \hfill $\Box$

Next we show that the equation (\ref{pde_0}) has a unique solution
and, consequently, $u_0^*=u_0|_{\Omega\setminus\overline{O}_0}$.
Assume that $u^1$ and $u^2$ are two solutions and let $u=u_2-u_1$,
then $u|_{\partial\Omega}=0$ and

\begin{eqnarray}0=&-\int_{\Omega\setminus\overline{O_\infty\cup O_0}}(\nabla
\cdot\sigma\nabla u)udx=-\int_{\partial\Omega}(\sigma \nabla u)\cdot \nu uds+\int_{\partial O_0} (\sigma \nabla u)\cdot \nu |_+uds\\&+\int_{\partial
O_\infty}(\sigma \nabla u)\cdot \nu |_+uds
+\int_{\Omega\setminus\overline{O_\infty\cup O_0}} |\nabla
u|_{\sigma}^2dx=\int_{\Omega\setminus\overline{O_\infty\cup O_0}} |\nabla u|_{\sigma}^2dx.
\end{eqnarray}
Thus $|\nabla u|\equiv 0$ in
$\Omega\setminus\overline{O}_0$. Since $\Omega\setminus\overline{O}_0$
is connected and $u=0$ at the boundary, we conclude uniqueness of
the solution of the equations (\ref{pde_0}). 

\begin{theorem}
Let $u_k$ and $u_0$ be the unique solution of (\ref{pde_k})
respectively (\ref{pde_0}) in $H^1(\Omega)$. Then
$u_k\rightharpoonup u$ and, consequently, $I_k[u_k]\xrightarrow{k\downarrow 0^+}I_0[u]$.
\end{theorem}

{\bf Proof:} We show first that $\{u_k\}$ is bounded in
$H^1(\Omega)$ uniformly in $k\in (0,1)$. Since $1/k>1$, we have
\begin{eqnarray*}
\frac{\lambda}{2}\|\nabla
u_k\|^2_{L^2(\Omega\setminus\overline{O}_0)}\leq\frac{1}{2}\int_{\Omega\setminus\overline{O_\infty \cup
O_0}} |\nabla u_k|_{\sigma}^2dx+\frac{1}{2k}\int_{O_\infty} |\nabla
u_k|_{\sigma_{1}}^2dx\\ \leq I_k[u_k]\leq I_k[u_0]\leq\frac{\Lambda}{2}\|\nabla
u\|^2_{L^2(\Omega\setminus\overline{O}_0)},
\end{eqnarray*}or
\begin{equation}\label{ineqb}
\|\nabla u_k\|^2_{L^2(\Omega\setminus
\overline{O}_0)}\leq\frac{\Lambda}{\lambda}\|\nabla
u\|^2_{L^2(\Omega\setminus\overline{O}_0)}.
\end{equation}
From (\ref{ineqb}) and the fact that $u_k|_{\partial\Omega}=f$, we
see that $\{u_k\}$ is uniformly bounded in $H^1(\Omega \backslash
O_0)$ and hence weakly compact. Therefore, there is a subsequence
$u_k\rightharpoonup u^*$ in $H^1(\Omega \backslash O_0)$, for some
$u^*$ with trace $f$ at ${\partial\Omega}$.

We will show next that $u^*$ satisfies the equations
(\ref{pde_0}), and therefore $u^*=u_0$ on $\Omega \setminus O_0$. By the
uniqueness of solutions of (\ref{pde_0}) we also conclude that the
whole sequence converges to $u$.

Since $u_k\rightharpoonup u^*$ we have that
%$0=\int_V\sigma_2\nabla u_k\cdot\nabla \phi dx\to \int_V\sigma_2\nabla u_0^*\cdot\nabla\phi dx$, for all $\phi\in C^\infty_0(V)$, which shows that
%$u_0^*$ solves $\nabla \cdot\sigma_2\nabla u_0^*=0$ in $V$. Similarly,
$0=\int_{\Omega\setminus \overline{O_0\cup O_\infty }}\sigma\nabla
u_k\cdot\nabla \phi dx\to \int_{\Omega\setminus \overline{O_\infty \cup
O_0}}\sigma \nabla u^*\cdot\nabla\phi dx$, for all $\phi\in
C^\infty_0({\Omega\setminus \overline{O_\infty\cup O_0}})$. Therefore $\nabla
\cdot\sigma\nabla u^*=0$ in ${\Omega\setminus \overline{O_\infty \cup
O_0}}$. Further, since $u_k$  minimizes $I[u_k]$ we must have
$\nabla u^*=0$ in $O_\infty$. To check the boundary conditions, note
that, for all $\phi\in C^\infty_0(\Omega)$ with $\phi\equiv 0$ in
$O_\infty$, we have $\int_{\partial O_0}(\sigma \nabla u_k)\cdot \nu |_+\phi ds=0$. Using the fact that $\phi$
were arbitrary, by taking the weak limit in $k\to 0$, we get
$\left.\frac{\partial u^*}{\partial\nu}\right|_+=0$ on $\partial
V$.  A similar argument applied to $\phi\in C^\infty_0(\Omega)$ with
$\phi\equiv 0$ in $O_0$, $\phi\equiv 1 $ in $O^j_{\infty}$, and $\phi\equiv 0$
in $O_\infty \backslash O_{\infty}^j$, also shows that $\int_{\partial
{O^j_{\infty}}}(\sigma \nabla u^*)\cdot \nu|_+\phi
ds=0$. Hence $u^*$ is the unique solution of the equation
(\ref{pde_0}) on $\Omega \backslash \overline{O}_0$. Thus $u_k$
converges weakly to the solution $u_0$ of (\ref{pde_0}) in
$\Omega\backslash \overline{O}_0$.

\section{Conclusions}
\label{conc}
We have considered the reconstruction of an anisotropic conductivity conformal  to a known $\sigma_0$ when one has knowledge of the internal functional $\sqrt{\sigma^{-1}_0 J \cdot J}$. Such data can be obtained by a novel combination of Current Density and Diffusion Tensor measurements. We have identified a variational problem defined in terms of the measured data and shown how to calculate the conformal factor from its unique solution. Further, we have presented a solution of the problem which allows for regions of infinite or zero conductivity, and which does not explicitly use Ohm's law in such regions.  We also proved that the equipotential sets  minimize the area functional corresponding to a Riemannian metric defined from the measured data.
\section{Acknowledgements}
The authors wish to kindly thank Mike Joy, Weijing Ma, and Nahla Elsaid for insightful discussions on DTI and the relation  between the conductivity and diffusion tensors.

\bibliographystyle{siam}

\end{document}